\documentclass{article}
\usepackage{amsmath,amsthm,amssymb}

\newcommand{\fc}{\mathfrak{c}}

\newcommand{\ga}{\alpha}
\newcommand{\gb}{\beta}

\newcommand{\gd}{\delta}
\newcommand{\gw}{\omega}

\newcommand{\gs}{\sigma}

\newcommand{\cP}{\mathcal{P}}

\newcommand{\dotrgen}{{\dot r}_{\mathit{gen}}}

\newcommand{\dom}{\mathrm{dom}}

\newcommand{\pioneoneonsigmaoneone}{$\mathbf{\Pi^1_1}$ on $\mathbf{\Sigma^1_1}$}
\newcommand{\pomodfin}{\cP(\gw)/\mathrm{fin}}
\newcommand{\gwgw}{[\gw]^{<\aleph_0}}

\newtheorem{theorem}{Theorem}[section]
\newtheorem{lemma}[theorem]{Lemma}
\newtheorem{claim}[theorem]{Claim}

\newtheorem{proposition}[theorem]{Proposition}

\theoremstyle{definition}
\newtheorem{definition}[theorem]{Definition}
\newtheorem{example}[theorem]{Example}
\newtheorem{question}[theorem]{Question}
\newtheorem{conjecture}[theorem]{Conjecture}
\newcommand{\trace}{\mathit{tr}}
\newcommand{\qi}{Q_{\trace(I)}}
\newcommand{\crn}{continuous reading of names}
\newcommand{\bintree}{2^{<\gw}}
\newcommand{\wtree}{\gw^{<\gw}}
\newcommand{\cantor}{2^\gw}
\newcommand{\baire}{\gw^\gw}
\newcommand{\exinfty}{\exists^\infty}
\newcommand{\liff}{\leftrightarrow}

\newcommand{\jx}{J^{<\gw}}
\newcommand{\dotagen}{\dot a_{gen}}

\title{Forcing with quotients}

\author{Michael Hru{\v s}{\' a}k
\thanks{Partially supported by GA {\v C}R grant 201-03-0933, PAPIIT grant IN108802-2 and CONACYT grant 40057-F.}\\
Universidad Nacional Aut{\' o}noma de M{\' e}xico\\
\and
Jind{\v r}ich Zapletal
\thanks{Partially supported by GA {\v C}R grant 201-03-0933, NSF grant DMS 0300201, and PAPIIT grant IN108802-2.
The results contained in the paper were obtained while the second author visited UNAM, Morelia, Mexico.}\\
University of Florida}

\begin{document}
\bibliographystyle{plain}
\maketitle
\begin{abstract}
We study an extensive connection between factor forcings of Borel subsets of Polish spaces modulo a $\gs$-ideal and factor
forcings of subsets of countable sets modulo an ideal.
\end{abstract}

\section{Introduction}

In recent years there has been a wave of interest in partial orders given as quotients, either Borel$(X)/I$ for a $\gs$-ideal $I$ on a Polish space $X$,
or Power$(X)/J$ for an ideal $J$ on a countable set $X$. The former turned out to be very close
to traditional forcings adding a real, and they allow of a comprehensive theory \cite{z:book}. The situation for the latter
is much less clear and less well understood. In this paper we describe a close relationship between the two classes of posets.
The connecting link is the following definition due to Brendle:

\begin{definition}
For a $\gs$-ideal $I$ on $\cantor$ or $\baire$ the \emph{trace ideal} $\trace(I)$ on $\bintree$ or $\wtree$ is defined by 
$a\in\trace(I)\liff\{r\in\cantor:\exinfty n\ r\restriction n\in a\}\in I.$
\end{definition}

We prove 

\begin{theorem}
\label{firsttheorem}
Suppose that $I$ is a $\gs$-ideal on $X=\cantor$ or $X=\baire$. If $P=$Borel$(X)/I$ is a proper forcing with continuous reading of names, 
then $Q=$Power$(\bintree$ or $\wtree)/\trace(I)$ is a proper forcing as well
and in fact $Q$ is naturally isomorphic to a two step iteration of $P$ and an $\aleph_0$-distributive
forcing.
\end{theorem}

This result makes it easy to generate and understand a large variety of factors of ideals on $\gw$. 
Our methods provide many ideals for which these factors are proper as well as examples of ideals for
which the factor forcings are improper.

Earlier results in this area include a note of Stepr{\= a}ns \cite{steprans:note}
on what in retrospect are trace ideals for a small class of forcings, a result of Balcar, Hern{\' a}ndez, and Hru{\v s}{\' a}k \cite{bhh:dense}
regarding the properness of the factor Power($\mathbb Q$)/nowhere dense sets, and results of Stepr{\= a}ns and Farah concerning the properness
of fators Power$(\gw)/J$ for various analytic P-ideals $J$. It should be noted that the trace ideals are analytic P-ideals only in the case
the original forcing $P_I$ had an exhaustive submeasure on it by a result of Solecki \cite{solecki:analytic}.

The notation used in the paper follows the set theoretic standard of \cite{jech:set}. If $I$ is a $\gs$-ideal on a Polish space $X$,
the symbol $P_I$ denotes the poset of $I$-positive Borel subsets of $X$ ordered by inclusion. If $J$ is an ideal on a countable set $X$, the symbol
$Q_J$ denotes the poset of $J$-positive subsets of $X$ ordered by inclusion. For a tree $T\subset(2\times\gw)^{<\gw}$ the symbol
$[T]$ denotes the set of all its infinite branches and the symbol $p[T]$ its projection, that is the set of those $r\in\cantor$ such that there is
$b\in\baire$ such that the pair $r,b$ constitutes a branch through the tree $T$. The characteristic function of a set $a\subset\gw$
is denoted by $\chi(a)$. For a sequence $t\in\bintree$ the symbol $[t]$ denotes the basic open subset of the space $\cantor$ determined by $t$.
LC denotes the use of suitable large cardinal assumptions.

\section{The \crn}

The following natural definition is critical for the wording and proofs of all results in this paper.

\begin{definition}
Let $I$ be a $\gs$-ideal on a Polish space. The forcing $P_I$ has the \emph{\crn}\ 
if for every $I$-positive Borel set $B$ and a Borel function $f:B\to\cantor$ there is an $I$-positive
Borel set $C\subset B$ such that the function $f$ is continuous on it.
\end{definition}

\noindent There are several natural restatements of this property.

\begin{claim}
\label{restatementclaim}
Let $I$ be a $\gs$-ideal on a Polish space $X$. The following are equivalent:

\begin{enumerate}
\item the forcing $P_I$ has the \crn 
\item for every $I$-positive Borel set $B$ and a countable collection $\{D_n:n\in\gw\}$
of Borel sets there is an $I$-positive Borel set $C\subset B$ such that all sets $D_n\cap C$ are relatively open in $C$
\item for every $I$-positive Borel set $B$ and every Borel function $f:B\to Y$ to a Polish space $Y$ there is an $I$-positive Borel set $C\subset B$
such that $f\restriction C$ is continuous.
\end{enumerate}
\end{claim}

\begin{proof}
(1)$\to$(2). Fix sets $B, D_n:n\in\gw$ and define a Borel function $f:B\to\cantor$ by $f(r)(n)=1$ if $r\in D_n$. By the
\crn\ there is an $I$-positive Borel set $C\subset B$ such that $f\restriction C$ is continuous. It is immediate that the sets
$D_n\cap C$ must be relatively open in $C$.

(2)$\to$(3). Suppose that $B$ is a Borel $I$-positive set and $f:B\to Y$ is a Borel function. For every basic open set
$O$ from some fixed countable basis for the space $Y$, let $D_O=f^{-1}O$. It is clear that $D_O$ is a Borel set and if $C\subset B$
is any set such that all sets $D_O\cap C$ are relatively open in $C$, the function $f\restriction C$ must be continuous.

(3)$\to$(1). Trivial.
\end{proof} 

\noindent Most definable proper partial orderings have the continuous reading of names under a suitable representation.

\begin{example}
\cite{z:book} 2.2.3. Every proper $\baire$-bounding poset $P_I$ has the \crn.
\end{example}

\begin{proof}
For simplicity assume that the underlying space of the ideal $I$ is $\cantor$. Suppose $B$ is a Borel $I$-positive set and $f:B\to\cantor$ is a Borel function.
Let $T\subset(2\times 2\times\gw)^{<\gw}$ be a tree which projects to the graph of the function $f$. By a standard absoluteness argument,
$B\Vdash$``for some $\dot s\in\cantor,
\dot b\in\baire$ the triple $\langle \dotrgen,\dot s, \dot b\rangle$ constitutes a branch through the tree $\check T$''. Since the forcing $P_I$ is bounding,
there is a condition $D\subset B$ which forces $\dot b$ to be pointwise dominated by some function $c\in\baire$. Let $S$ be the subset of the tree $T$
consisting of those sequences whose third coordinate is pointwise dominated by the function $c$. Then $S$ is a finitely branching tree and

\begin{itemize}
\item $p[S]$ is a compact subset of the graph of the function $f$, so it is a graph of a continuous subfunction of $f$
\item $C=pp[S]$ is a compact subset of the set $B$, $D$ forces the generic real into $\dot C$ and therefore $C$ is $I$-positive.
\end{itemize}

All in all, $C\subset B$ is an $I$-positive compact set on which the function $f$ is continuous.
\end{proof}

\begin{example}
If the ideal $I$ is $\gs$-generated by closed sets then $P_I$ has the \crn.
\end{example}

\begin{proof}
Let us first deal with the case of the $\gs$-ideal generated by nowhere dense sets. It is a classical result \cite{kechris:book} 8.38 that for every Polish space
$X$ and every Borel function $f:X\to\baire$ there is a comeager $G_\gd$ set $C\subset X$ such that $f$ is continuous on it.

In the general case, suppose that $B$ is a Borel $I$-positive set and $f:B\to\baire$ is a Borel function. By a result of Solecki
\cite{solecki:analytic}, thinning out the set $B$ we may assume that it is $G_\gd$ and for every basic open set $O$, $O\cap B\notin I\liff O\cap B\neq 0$.
It immediately follows that then every closed set in the ideal $I$ is nowhere dense in the set $B$. Since the set $B$ is $G_\gd$, it is Polish
in the relative topology, and so every set $C\subset B$ comeager in $B$ must be positive in the ideal $I$. By the first paragraph of the proof,
there must be a comeager $G_\gd$ set $C\subset B$ such that the function $f$ is continuous on it.
\end{proof} 

These two classes of examples include many forcings used in practice, such as the Cohen, Solovay, or Miller reals. In other situations, the \crn\ 
has to be checked carefully.

\begin{example}
The Laver forcing in the natural presentation has the \crn.
\end{example}

\begin{example}
\label{stepransexample}
The \emph{Stepr\=ans forcing} \cite{steprans:discontinuous} in the natural presentation does not have the \crn. Here 
the Stepr{\= a}ns forcing $P_I$ is obtained from a Borel function $f:\cantor\to\cantor$
which cannot be decomposed into countably many continuous functions by considering the ideal $I$ $\gs$-generated by
the sets on which the function $f$ is continuous. The poset $P_I$ is proper and up to the forcing equivalence does not depend on the initial choice of
the function $f$--see \cite{z:book}, 2.3.49.

It is interesting to note that in a slightly different presentation the ideal associated with Stepr{\= a}ns forcing is generated by closed sets and therefore does
have the \crn. Let us describe this different presentation. 

We will need a definite example of a Borel function which cannot be decomposed into countably many continuous functions, due to Pawlikowski \cite{cmps:decompose}.
Consider the space $\gw+1$ equipped with the order topology, the space $(X,\gs)=(\gw+1)^\gw$ with the product topology, and the \emph{Pawlikowski function}
$P:X\to\baire$ defined by $P(r)(n)=r(n)+1$ if $r(n)\in\gw$ and $P(r)(n)=0$ if $r(n)=\gw$. This function cannot be decomposed into countably many
continuous functions and it is in a sense a minimal such example \cite{solecki:continuous}.

Let $I$ be the $\gs$-ideal on the space $X$ generated by the sets $B$ on which the function $P$ is continuous. Clearly the poset $P_I$ does not
have the \crn\ as witnessed by the function $P$. However, the function $P$ turns out to be the only obstacle. Namely, if the space $X$ is equipped with the smallest
Polish topology $\tau\supset\sigma$ which makes the function $P$ continuous and generates the same Borel structure, the $\gs$-ideal $I$ is generated
by $\tau$-closed sets and so the forcing $P_I$ has the \crn\ in this new presentation. An outline of the easy argument: the topology $\tau$ is the
product topology on $X$ if $\gw+1$ is viewed as a discrete space. If $B\subset X$ is a set such that $P\restriction B$ is continuous with respect to
the topology $\gs$ and $C$ is the $\tau$-closure of the set $B$, then $P\restriction C$ is continuous with respect to $\gs$ as well. If $U,V$ are basic open subsets
of $(X,\gs)$ and $\baire$ respectively such that $P''(B\cap U)\subset V$ then $P''(C\cap U)\subset V$ as well.
\end{example}

\begin{definition}
Let  $J$ be an ideal on $\gw$. The \emph{Prikry forcing} $P(J)$ for the ideal $J$ is defined as the set of all pairs
$\langle t, a\rangle$ where $t\subset\gw$ is a finite set, $a\subset\gw$ is a set in the ideal $J$, and $\langle u,b\rangle\leq \langle t,a\rangle$
if $t\subset u, a\subset b$ and $a\cap u\setminus t=0$. We will refer to the union of the first coordinates of conditions in the generic filter as the
generic subset of $\gw$, and denote it by $\dotagen$.
\end{definition}

\begin{example}
\label{prikryexample}
Let $J$ be an ideal on $\gw$. The forcing $P(J)$ has the \crn\ if and only if $J$ is a P-ideal.
\end{example}

\begin{proof}
Let $I$ be the $\gs$-ideal on $\cantor$ associated with the forcing $P(J)$, namely $I$ is the collection of those sets $B\subset\cantor$ for which
it is outright forced that $\chi(\dotagen)\notin B$. Thus the poset $P_I$ is in the forcing sense equivalent
to the poset $P(J)$, with a canonical correspondence between the respective generic objects.

First suppose that $J$ is not a P-ideal, as witnessed by a countable collection $\{a_n:n\in\gw\}$ of sets in the ideal such that no
set in the ideal contains each of them modulo a finite set. Consider
the Borel function $f:\cantor\to\cantor$ defined by
$f(r)(n)=$the parity of the size of the set $\{m\in a_n:g(m)=1\}$. The function $f$ is defined on an $I$-large set,
and we claim that it cannot be reduced to a continuous function on an
$I$-positive Borel set. 

Suppose that $B$ is an $I$-positive Borel set, and $\langle t,b\rangle\Vdash$``$\chi(\dotagen)\in\dot B$''. 
Thinning out the set $B$ we may assume that it consists only of functions $r$ such that $\forall m\in\gw\ (t(m)=1\to r(m)=1$ and $m\in b\to r(m)=0)$.
Let $n$ be such that the set $a_n\setminus b$ is infinite.
It is not difficult to see that both sets $\{r\in B:f(r)(n)=0\}$ and $\{r\in B:f(r)(n)=1\}$ are dense in the set $B$,
and therefore the function $f$ cannot be continuous on $B$.

Now suppose that $J$ is a $P$-ideal, $B\notin I$ is a Borel set, and $f:B\to\baire$ is a Borel function. Let $M$ be a countable elementary submodel
of a large enough structure, let $a\subset\gw$ be a set in the ideal $J$ which modulo finite
contains all sets in $J\cap M$, and for every number $n$ consider the sets $C_n=\{r\in B:r$ is $M$-generic and for all $k>n$, $r(k)=1\to k\notin a\}$.
Since the poset $P_I$ is c.c.c., the set $B\setminus\bigcup_nC_n$ is in the ideal $I$ and there must be a number $n\in\gw$ such that
the Borel set $C_n$ is $I$-positive. Set $C=C_n\subset B$; we will be done if we show that
the function $f\restriction C$ is continuous.

Suppose $r\in C$ and $O\subset\baire$ is a basic open set such that $f(r)\in O$. We must produce a basic open set $P\subset\cantor$ such that
$r\in P$ and for every real $s\in P\cap C$, $f(r)\in O$. Look at the $M$-generic filter $G\subset M\cap P(J)$ associated with the real $r$:
$G=\{\langle t,b\rangle\in P(J)\cap M: \forall m\in\gw\ m\in t\to r(m)=1\land m\in b\to r(m)=0\}$. By the forcing theorem, there must be a condition $\langle t,b\rangle\in G$
such that it forces $\dot f(\chi(\dotagen))\in O$. Let $m\in\gw$ be a natural number larger than $n$,
larger than all elements of the finite set $t$, and larger than all elements of the finite set $b\setminus c$. It is enough to show that
whenever $s\in C$ is a real such that $s\restriction m=r\restriction m$ then $f(s)\in O$. A brief inspection reveals that the condition
$\langle t,b\rangle$ belongs to the $M$-generic ultrafilter associated with the real $s$, and by the forcing theorem applied in the model $M$,
it must be the case that $f(s)\in O$ as desired. 
\end{proof}

A similar proof can be used to show that the Hechler forcing in the natural presentation has the \crn,
while the eventually different real forcing does not have the \crn.

\begin{example}
The \emph{eventually different real forcing} \cite{bartoszynski:set} 7.4.8
does not have the \crn. The forcing is the set of all pairs $\langle t,x\rangle$ where $t\in\wtree$ and $x\subset\baire$
is a finite set, and $\langle u, y\rangle\subset\langle t,x\rangle$ if $t\subset u$, $x\subset y$, and $\bigcup x\cap u\setminus t=0$. The generic function
is the union of the first coordinates of the conditions of the generic filter, and it is eventually different from all functions in $\baire$ in the ground model.
Let $I$ be the $\gs$-ideal associated with this forcing. As in the previous proof, the function $f(r)(n)=$the parity of the size of the set $\{m\in\gw:r(m)=n\}$
is a function which cannot be reduced to a continuous function on an $I$-positive Borel set.
\end{example}

The \crn\ is a rather slippery property of ideals. It is not preserved under Borel isomorphism of ideals. This is to say that
there are $\gs$-ideals $I$ and $J$ on Polish spaces $X$ and $Y$ and a Borel bijection $f:X\to Y$ such that a set $A\subset Y$
is in the ideal $J$ iff its $f$-preimage is in the ideal $I$, but the poset $P_I$ does have the \crn\ while $P_J$ does not.
An instructive example is that of the Stepr{\= a}ns forcing,~\ref{stepransexample}.
Note that since Borel injective images of Borel sets are Borel, in this case the function $f$
can be naturally extended to an isomorphism of the posets $P_I$ and $P_J$. This means that the \crn\ is, in fact, a property of a presentation of forcing as opposed
to a property of the forcing itself. Even so, the \crn\ is perceived as a natural and useful property. We state two of its important features.

\begin{claim}
(LC) If $I, J$ are universally Baire $\gs$-ideals on Polish spaces such that the forcings $P_I,P_J$ are provably 
proper and have the continuous reading of names then the iteration
$P_I*P_J$ in the natural presentation has the continuous reading of names.
\end{claim}

\begin{proof}
For simplicity assume that the Polish spaces in question are both just $\cantor$. To obtain the natural presentation of the iteration,
use the work of \cite{z:book}. Add all the universally Baire sets without Borel $I$-positive (or $J$-positive, respectively) subsets to the ideal $I$ or $J$ respectively.
Consider the Fubini product $I*J$ of the two ideals on the space $\cantor\times\cantor$. Then $P_I*P_J$ is in the forcing sense naturally
isomorphic to $P_{I*J}$. Moreover, every $I*J$-positive Borel set contains a \emph{good} Borel subset, that is a set $B$ such that 
its projection $p(B)$ into the first coordinate is Borel and for every real $r\in p(B)$ the vertical section $B_r$ is $J$-positive.

So let $B\subset \cantor\times\cantor$ be a good $I*J$-positive Borel set and $f:B\to\cantor$ be a Borel function. For every real $r\in p(B)$
there is a $J$-positive Borel subset of the vertical section $B_r$ on which the function $f$ is continuous. Using universally Baire selection
and \cite{z:book} 2.2.1(2), there is a good $I*J$-positive set $C\subset B$ such that for every real $r\in p(C)$, the function $f$
is continuous on the vertical section $C_r$ and moreover the function $g:p(C)\to\cP(\bintree\times\bintree)$ is Borel, where $\langle u,v\rangle\in g(r)$
if for every real $s\in C_r$ extending $u$ the functional value $f(r,s)$ extends $v$. Let $D\subset p(C)$ be a Borel $I$-positive set
on which the function $g$ is continuous. Let $E\subset C$ be the set $\{\langle r,s\rangle\in C:r\in D\}$. It is not difficult to verify that
the set $E$ is Borel, $I*J$-positive, and the function $f\restriction E$ is continuous.
\end{proof}

We do not know if the \crn\ is preserved under (countably transfinite) (countable support) iterations of (definable) proper forcing notions, even though
in all natural particular cases this can be manually checked to hold true.

\begin{claim}
\cite{z:book} 2.2.2(2) Suppose that $I$ is a $\gs$-ideal on a Polish space. If $P_I$ is a proper forcing notion with the \crn\ 
then every Borel $I$-positive set has a $G_\gd$ $I$-positive subset.
\end{claim}

\begin{proof}
For simplicity assume that the underlying space is $\cantor$. Suppose that $B$ is an $I$-positive Borel set and $T$ is a tree on $2\times\gw$
projecting to it. Since the poset $P_I$ is proper, there must be an $I$-positive Borel set $D\subset B$ and a Borel function $f:D\to[T]$
such that for every element $r\in D$, the first coordinate of the branch $f(r)$ of the tree is $r$ itself \cite{z:book} Lemma 2.2.1(2). By the \crn\ we may assume that
the function $f$ is continuous on $D$. Every partial continuous function can be extended to a continuous function with
a $G_\gd$ domain. Let $D\subset C, f\subset g$ be such a $G_\gd$ set and a continuous extension, with $D$ still dense in $C$.
It is immediate that $g:C\to[T]$ is a function such that for every $r\in C$ the first coordinate of the branch $g(r)$ of the tree $T$ is
$r$ itself. Then $C\subset B$ is an $I$-positive $G_\gd$ subset of the set $B$.
\end{proof}

The opposite implication does not hold: compact sets are dense in the natural presentation of Stepr{\= a}ns forcing \cite{z:book} 2.3.46, while the \crn\ fails.

\section{Proof of Theorem~\ref{firsttheorem}}

We will deal with a $\gs$-ideal $I$ on the space $\cantor$, the $\baire$ case being analogous.

\begin{definition}
The function $\pi:\cP(\bintree)\to\cP(\cantor)$ is defined by $\pi(a)=\{r\in\cantor:\exinfty n\ r\restriction n\in a\}$.
\end{definition}

Clearly, the range of the function $\pi$ is exactly the collection of all $G_\gd$ subsets of $\cantor$, and the function $\pi$ preserves inclusion.
Moreover, if $I$ is a $\gs$-ideal on $\cantor$ then $a\in\trace(I)$ if and only if $\pi(a)\in I$, and the map $\pi\restriction \qi:\qi\to P_I$ preserves compatibility.
For the remainder of the section fix a $\gs$-ideal $I$ on $\cantor$ such that the poset $P_I$ is proper and has the \crn, and write $J=\trace(I)$.

\begin{claim}
\label{embeddingclaim}
\begin{enumerate}
\item For every set $a\notin Q_J$ and for every $I$-positive $G_\gd$ subset $B\subset\pi(a)$ there is a set
$b\subset a$ such that $\pi(b)=B$.
\item $Q_J$ forces $\pi''\dot G$ to be a $P_I$-generic filter, where $\dot G$ is the name for the $Q_J$ generic filter.
\end{enumerate}
\end{claim}

\begin{proof}
The second item immediately follows from the first.
For the first one, suppose $B\subset \pi(a)$ is an $I$-positive $G_\gd$ set, $B=\bigcap_nO_n$ for some open sets $O_n$. By induction on $n\in\gw$ build sets
$a_n\subset a$ in the following way:

\begin{enumerate}
\item Each $a_n\subset\bintree$ is an antichain and it refines $a_{n-1}$. For notational convenience let $a_{-1}=\{0\}$.
\item $B\subset\bigcup_{t\in a_n}[t]\subset O_n$.
\end{enumerate}

After the construction is complete, clearly writing $b=\bigcup_na_n$ we will have $\pi(b)=B$ as required.
Suppose the antichain $a_n$ has been obtained. For each $t\in a_n$ let $c(t)$ be the collection of all proper extensions $u\in\bintree$ of $t$ such that the basic open set
determined by $u$ is a subset of $O_{n+1}$, and no proper initial segment of $u$ longer than $t$ has this property. 
Note that $c(t)\subset\bintree$ is an antichain. For each $u\in c(t)$ let $d(u)$ be the collection of all proper extensions $v$ of $u$ which are in the set $a$,
such that no proper initial segment of $v$ longer than $u$ has this property. Note that each $d(u)$ is an antichain. It is not difficult to verify
that the set $a_{n+1}=\bigcup_u d(u)$ has the desired properties.
\end{proof}

\begin{claim}
\label{propernessclaim}
The poset $Q_J$ is proper.
\end{claim}

\begin{proof}
Let $M$ be a countable elementary submodel of a large structure and let $a\in Q_J\cap M$ be a condition. Let $\langle D_n:n\in\gw\rangle$ be an enumeration of all
open dense subsets of the poset $Q_J$ in the model $M$. We will find sets $a_n\subset a$ and functions $g_n:a_n\to D_n\cap M$ so that

\begin{enumerate}
\item each set $a_n\subset\bintree$ is an antichain and it refines $a_{n-1}$. For notational convenience let $a_{-1}=\{0\}$.
\item the set $b=\bigcup_na_n\subset a$ is $J$-positive
\item for each sequence $t\in a_n$ the set $\{u\in b:t\subset u\}$ is a subset of $g_n(t)$.
\end{enumerate}

It immediately follows that the set $b\subset a$ is the required $M$-master condition in the poset $Q_J$. To see this, choose a $J$-positive set $c\subset b$
and an open dense set $D=D_n\in M$ for some number $n$. For each sequence $t\in a_n$ write $b_t=\{u\in b:t\subset u\}$. Since the set
$a_n\subset\bintree$ is an antichain, it is the case that $\pi(c)=\bigcup_{t\in a_n}\pi(c\cap b_t)$ and therefore one of the sets $c\cap b_t:t\in a_n$
must be $J$-positive. Such set $c\cap b_t\subset c$ has the condition $g_n(t)\in D_n\cap M$ above it as required.

To perform the construction, find an $M$-master condition $B\subset\pi(a)$ for the poset $P_I$. Thinning out the condition $B$ we may assume that for
every dense set $E\in M$ of the poset $P_I$, $B\subset\bigcup(E\cap M)$. Thinning out the condition $B$ even further, by Claim~\ref{restatementclaim},
we may assume that for every set $C\in P_I\cap M$ the intersection $C\cap B$ is relatively open in $B$. Thinning out the condition $B$ further
we may assume that it is a $G_\gd$ set such that for every basic open set $O$, $B\cap O\notin I\liff B\cap O\neq 0$.
Fix a representation $B=\bigcap_n O_n$, for some open sets $O_n$.

The induction hypotheses on the construction of the sets $a_n$ are the following.

\begin{enumerate}
\item Each $a_n\subset\bintree$ is an antichain and it refines $a_{n-1}$. For notational convenience let $a_{-1}=\{0\}$.
\item $B\subset\bigcup_{t\in a_n}[t]\subset O_n$
\item For every $n$, $g_n(t)\subset \{u\in a:t\subset u\}$ is a condition in the open dense set $D_n$.
For every $n\in m$, $t\in a_n$ and $u\in a_m$ such that $t\subset u$, $u\in g_n(t)$ and $g_m(u)\subset g_n(t)$. For notational convenience
put $g_{-1}(0)=a$.
\item For each $t\in a_n$ it is the case $B\cap[t]$ is a nonempty subset of $\pi(g_n(t))$.
\end{enumerate}

Now suppose that $a_n, g_n$ have been constructed. Fix a node $t\in a_n$. We will show how the part of the antichain $a_{n+1}$ below $t$ will be constructed.
Let $E$ be the part of the open dense set $D_{n+1}\subset Q_J$ below the condition $g_n(t)\in Q_J$. Claim~\ref{embeddingclaim} shows that the set
$\pi''(E)$ is dense below the condition $\pi(g_n(t))$. Then $B\cap[t]\subset\bigcup(\pi''E\cap M)=\bigcup\pi''(E\cap M)$
by the choice of the $M$-master condition $B$. Note that for every condition $p\in E\cap M$ the set $\pi(p)\cap B\subset B\cap[t]$ is relatively open
by the choice of the condition $B\in P_I$ again.
It is now easy to build an antichain $b\subset\bintree$ below the node $t$ so that for every $u\in b$ it is the case that $[u]\subset O_n$ and there is
a condition $p(u)\in E\cap M$ such that $B\cap[u]$ is a nonempty subset of $\pi(p(u))$, and $B\cap[t]\subset\bigcup_{u\in b}[u]$. Let
$c$ be then the collection of all nodes $v\in\bintree$ such that there is some $u\in b$ such that $u\subset v$, $v\in p(u)$, $B\cap[v]\neq 0$ and no
proper initial segment of $v$ is an extension of $u$ in $p(u)$. The set $c$ is an antichain below the node $t$, and it is the part of the
antichain $a_{n+1}$ below $t$. For every node $v\in c$ let $g_{n+1}(v)=p(u)\cap\{w\in\bintree:v\subset w\}$. The induction hypotheses are
easily seen to be satisfied.
\end{proof}

\begin{claim}
The remainder poset $R=Q_J/P_I$ preserves stationary subsets of $\gw_1$ and it is $\aleph_0$-distributive.
\end{claim}

\begin{proof}
Here the remainder poset $R$ is computed via the $Q_J$-name for the $P_I$-generic filter obtained in Claim~\ref{embeddingclaim}.
Note that writing $\dotrgen$ for the canonical $P_I$-generic real we have $t\subset\dotrgen\liff$ the set $\{u\in\bintree:t\subset u\}$
is in the $Q_J$-generic filter, this for every sequence
$t\in\bintree$. 

The fact that $P_I\Vdash$``$\dot R$ is stationary preserving'' follows abstractly from the proof of the previous claim: if $M$ is a countable elementary submodel
and $B\in P_I$ is any $M$-master condition for the poset $P_I$ then there is an $M$-master condition $b\in Q_J$ such that $\pi(b)\subset B$.
For the distributivity, suppose that $\dot f$ is a $Q_J$-name for an $\gw$-sequence of ordinals. We must prove that $\dot f\in V[\dotrgen]$.

To this end, revisit the proof of the previous claim again. Assume that $\dot f\in M$ and for each number $k\in\gw$ find a number $n_k\in\gw$
such that the conditions in the open dense set $D_{n_k}\subset Q_J$ decide the value of $\dot f(\check k)$. Look again at the master
condition $b=\bigcup_na_n$. It is not difficult to see that $b$ forces that for each $n\in\gw$ there is exactly one initial segment of the real $\dotrgen$
in the set $a_n$; call it $t_n$. Consequently, the sequence $\dot f$ can be recovered in the model $V[\dotrgen]$ by the following formula:
$\dot f(\check k)$ is that ordinal which is forced by the condition $g_{n_k}(t_{n_k})$ to be the value of $\dot f(\check k)$.
\end{proof}

It is interesting to see what the $\aleph_0$-distributive tail $Q_J/P_I$ can be. From the definitions it is
equal to the collection of all ground model sets $a\subset\bintree$ such that the $P_I$-generic real has infinitely many initial segments
in $a$, ordered by inclusion. In many cases it is, in the forcing sense, equivalent to 
$\pomodfin$ of the $P_I$ extension. To prove this it is just enough to show that $P_I$ forces every infinite subset of the generic real $\dotrgen$ 
(understood now as a path through $\bintree$
or $\wtree$) to have an infinite subset of the form $a\cap\dotrgen$ for some set $a$ in the ground model. We can verify this property in a great number
of cases and disprove in others, but we do not have a suitable general criterion.

\begin{proposition}
Let $I$ be a $\gs$-ideal on $\cantor$ $\gs$-generated by a $\gs$-compact family of closed sets. The forcing $P_I$ is proper and bounding, and
writing $J=\trace(I)$, $Q_J=P_I*\pomodfin$.
\end{proposition}

\begin{proof}
The ideals $I$ considered in this proposition form a class considered in \cite{z:four}. There it is proved that the poset $P_I$ is proper and bounding.
In fact a standard determinacy argument \cite{z:four} shows the following. Fix a $\gs$-ideal $I$ $\gs$-generated by a $\gs$-compact collection $F=\bigcup_n F_n$
of closed sets, where the sets $F_n\subset K(\cantor)$ are closed. Call a tree $T\subset\bintree$ $I$-\emph{good} if for every node $t\in T$ and every
number $n$ there is a number $m$ such that no set in $F_n$ meets all the open sets determined by the extensions of the node $t$ in the tree $T$ of length $m$.
Then a Borel set $B\subset\cantor$ is $I$-positive if and only if it contains all branches of some $I$-good tree. Therefore the poset of $I$-good
trees is naturally isomorphic with a dense subset of the poset $P_I$ and below we will identify it with $P_I$.

We will show that $P_I\Vdash$``every infinite subset of $\dotrgen$ has an infinite subset of the form $a\cap\dotrgen$ for some set $a$ in the ground model''.
Suppose $T\in P_I$ is an $I$-good tree, $T\Vdash$``$\dot x\subset\dotrgen$ is an infinite set''. A standard fusion argument will give 
an $I$-good tree $S\subset T$ such that for every number $n$ there is $m>n$
such that every sequence $s\in S$ of length $m$ has an initial segment of length $\geq n$
in the set $a=\{t\in S:S\restriction t\Vdash\check t\in\dot x\}\subset S$. Clearly
$S\Vdash$``$\check a\cap\dotrgen$ is an infinite subset of $\dot x$'' as desired.
\end{proof}

The following definition is not standard. It is an attempt to restate a commonly used combinatorial forcing property in topological terms.

\begin{definition}
Let $I$ be a $\gs$-ideal on some Polish space $X$ with a fixed metric $d$. We say that the poset $P_I$ has the \emph{pure decision property} (with respect to
the metric $d$) if for every $I$-positive Borel set $B\subset X$ and every Borel map $f:(B,d)\to(Y,e)$ into a compact metric space
there is a Borel $I$-positive set $C\subset B$ on which the map $f$ is a contraction.
\end{definition}

\begin{example}
The Laver forcing has the pure decision property in the standard representation, with respect to the metric of least difference on $\baire$: $d(x,y)=2^{-n}$
where $n$ is the smallest number where the functions $x,y\in\baire$ differ.
\end{example}

\begin{proof}
Let $B$ be Borel $I$-positive set and $f:(B,d)\to (Y,e)$ be a Borel map into a compact metric space. Thinning out the set $B$ if necessary
we may assume that $B=[T]$ for some Laver tree $T\subset\wtree$. To simplify the notation assume that $T$ has an empty trunk.

Before we proceed recall the well known fact that for every Laver tree $S$ and Borel partition $[S]=\bigcup_{i\in n}A_i$ into finitely many pieces
there is a Laver tree $U\subset S$ with the same trunk such that the set $[U]$ is included in one of the pieces of the partition.

Now for every $n$ find a finite $2^{-n-1}$-network $y_n\subset Y$, that is, a set such that
every point of the space $Y$ is $2^{-n-1}$-close to one of its elements. By induction on $n\in\gw$ build a fusion sequence of Laver trees $T_n$
so that $T_0=T,$ $T_{n+1}$ agrees with $T_n$ on sequences of length $n+1$ and for every such a sequence $t\in T_n$ there is an element
$x_t\in y_n$ such that for every path $r$ through $T_{n+1}$ extending the sequence $t$, the element $f(r)\in Y$ is $2^{-n-1}$-close to $x_t$.
This is possible by the observation in the previous paragraph. Note that by the triangle inequality this means that for two such paths
$r_0, r_1$ the elements $f(r_0), f(r_1)\in Y$ will have $e$-distance $\leq 2^{-n}$. Let $S$ be the fusion of the sequence of trees $T_n$.
It is not difficult to see that the set $C=[S]$ has the required properties.
\end{proof}

\begin{proposition}
If $I$ is a $\gs$-ideal on $\baire$ such that the poset $P_I$ is proper and has the pure decision property
with respect to the metric of least difference on $\baire$. Then $\qi=P_I*\pomodfin$.
\end{proposition}

\begin{proof}
Note that the pure decision property implies the \crn.

Suppose that $B\in P_I$ forces $\dot x\subset\dotrgen$ to be an infinite set. Since the poset $P_I$ is proper, thinning out the set $B$ if necessary we can
find a Borel map $f:B\to\cantor$ such that $B\Vdash$``$\dot x=\{\dotrgen\restriction n:n\in \dot f(\dotrgen)\}$''. Consider the metric $e$ of least difference 
on $\cantor$ and
use the pure decision property to find an $I$-positive set $C\subset B$ such that $f:C\to\cantor$ is a contraction. This means that for every sequence $t\in\wtree$,
all reals $r\in C$ extending the sequence $t$ return the same value $b(t)\in 2$ for $f(r)(|t|)$. Let $a=\{t\in\wtree:b(t)=1\}$. It follows from the definitions
that $C\Vdash\dot x=a\cap\dotrgen$, and the lemma follows.
\end{proof}

\begin{example}
The Cohen poset forces that there is an infinite set $x\subset\dotrgen$ without an infinite subset of the form $a\cap\dotrgen$, $a\in V$. Just let
an initial segment $t$ of $\dotrgen$ into $\dot x$ if and only if $\dotrgen(|t|)=0$.
\end{example}

As a final remark in this section, once we produced so many ideals $J$ for which the factor forcing $Q_J$ is proper, we should also produce
some for which it is not proper. The following Proposition of independent interest shows how to do exactly that in several ways. First, an instrumental definition.

\begin{definition}
Let $\gb$ be a limit ordinal. We say that an inclusion-decreasing sequence $\langle I_\ga:\ga\in\gb\rangle$  of $\gs$-ideals on a Polish space \emph{does not stabilize}
if for every ordinal $\ga\in\gb$ and every $I_\ga$-positive Borel set $B$ there is an ordinal $\ga\in\gamma\in\gb$ and a Borel set $C\subset B$ which is $I_\ga$-small but
$I_\gamma$-positive. This is equivalent to saying that, writing $I=\bigcap_\ga I_\ga$, the sets $I_\ga\cap P_I$ are all dense in $P_I$.
Restated again, this is equivalent to saying that for every $I$-positive Borel set $B$, $I\restriction B\neq I_\ga\restriction B$--hence the terminology.

Similarly, we say that an inclusion-decreasing sequence $\langle J_\ga:\ga\in\gb\rangle$ of ideals on some countable set $X$ does not stabilize if for every
ordinal $\ga\in\gb$ and $J_\ga$-positive set $a\subset X$ there is an ordinal $\ga\in\gamma\in\gb$ and a set $b\subset a$ which is $J_\ga$-small but
$J_\gamma$-positive. This is the same as to say, writing $J=\bigcap_\ga J_\ga$, that the set $J_\ga\cap Q_J$ are dense in the factor forcing $Q_J$.
\end{definition} 

\begin{proposition}
Assume the Continuum Hypothesis. If $I$ is a $\gs$-ideal on a Polish space, then

\begin{enumerate}
\item $P_I$ collapses $\aleph_1$ if and only if $I=\bigcap_{n\in\gw}I_n$ for an inclusion-decreasing sequence of $\gs$-ideals which does not stabilize.
\item Suppose $P_I$ preserves $\aleph_1$. $P_I$ is nowhere c.c.c. if and only if $I=\bigcap_{\ga\in\gw_1}I_\ga$ for an inclusion-decreasing sequence
of $\gs$-ideals which does not stabilize.
\end{enumerate}

If $J$ is an ideal on a countable set, then

\begin{enumerate}
\item[3.] $Q_J$ adds an unbounded real if and only if $J=\bigcap_{n\in\gw}J_n$ for an inclusion-decreasing sequence of ideals which does not stabilize.
\item[4.] If $J=\bigcap_{n\in\gw}J_n$ for an inclusion-decreasing sequence of P-ideals which does not stabilize, then $Q_J$ collapses $\aleph_1$.
\end{enumerate}
\end{proposition}

\begin{proof}

For the first equivalence, assume that $P_I$ collapses $\aleph_1$. Let $\dot f:\check\gw\to\check\gw_1$ be a name for a function with cofinal range. 
For every number $n\in\gw$ let $I_n$ be the ideal generated by sets $B\in P_I$ which force the first $n$ values of the function $\dot f$ to be bounded by some fixed
countable ordinal, together with all sets in the ideal $I$.
It is clear that $\langle I_n:n\in\gw\rangle$ is an inclusion-decreasing sequence of $\gs$-ideals which does not stabilize, and $I=\bigcap_nI_n$.
On the other hand, suppose that $I=\bigcap_nI_n$ for some inclusion-decreasing sequence of $\gs$-ideals which does not stabilize. Since the ideals $I_n$ are dense in the
poset $P_I$, we can pick a maximal antichain $A_n\subset I_n$ from each, and by CH it will be enough to show that every condition in $P_I$
is compatible with uncountably many elements of one of these antichains. Indeed, if $B\in P_I$ is a condition, then $B\notin I_n$ for some number $n$,
and $B$ must be compatible with uncountably many elements of the antichain $A_n$, because if $X\subset A_n$ is a countable set, then $C=\bigcup X\in I_n$ and the
condition $B\setminus C\notin I_n$ is a condition incompatible with all elements of the set $X$.

For the second equivalence, assume that the poset $P_I$ preserves
$\aleph_1$ and is nowhere c.c.c. Then there is a name $\dot f$ for a function from $\gw_1$ to itself which is not bounded
by any ground model such function. To see this, let $\langle M_\ga:\ga\in\gw_1\rangle$ be a tower of countable elementary submodels of some large structure, and define
$\dot f(\ga)=\min\{\gb\in\gw_1:$ for every maximal antichain $A\in M_\ga$ the unique element in it which belongs to the generic filter is in the model $M_\gb\}$.
Since the forcing $P_I$ preserves $\aleph_1$, and by CH $P_I\subset\bigcup_\ga M_\ga$, this is well-defined. If $p\in P_I$ is a condition and
$g:\gw_1\to\gw_1$ is a function, find an ordinal $\ga\in\gw_1$ such that $p\in M_\ga$, a maximal antichain $A\in M_\ga$ which has uncountably
many elements below the condition $p$, and an element $q\in A\setminus M_{g(\ga)}$ below the condition $p$. Then $q\Vdash$``$g(\ga)\in\dot f(\ga)$'', and
it follows that the function $\dot f$ is unbounded.

Now, given an ordinal $\ga$ let $I_\ga$ be the ideal generated by the sets $B\in P_I$ for which 
there is a countable ordinal $\gb$ such that $B$ force all values $\{\dot f(\gamma):\gamma\in\ga\}$ to be smaller
than $\gb$, together with all sets in the ideal $I$. 
It is clear that $\langle I_\ga:\ga\in\gw_1\rangle$ is an inclusion decreasing sequence of $\gs$-ideals. 
Since the function $\dot f$ is not dominated by any ground model function, it is the case that $I=\bigcap_\ga I_\ga$, and since the
forcing $P_I$ preserves $\aleph_1$, the sequence of ideals does not stabilize.

For the other direction, let $I=\bigcap_\ga I_\ga$. Suppose $B\in P_I$ is a Borel set; we must find an uncountable antichain below it. It must be the case
that $B\notin I_\ga$ for some countable ordinal $\ga$. Now since the sequence of ideals does not stabilize, the ideal $I_\ga$ is dense in the poset $P_I$,
and therefore there must be a maximal antichain $A$ below $B$ which consists solely of $I_\ga$-small sets. This antichain must be uncountable, because otherwise
$\bigcup A\in I_\ga$ and $B\setminus\bigcup A\notin I_\ga$ is a condition in $P_I$ which avoids all elements of the maximal antichain, a contradiction.

For the third equivalence, first suppose that $Q_J\Vdash$``$\dot f\in\baire$ is an unbounded function''. Let $J_n=\{a\subset\gw:$ there is
a number $m$ forces the first $n$ values
of the function $\dot f$ to be smaller than $m\}$. It is immediate that $\langle J_n:n\in\gw\rangle$ is an inclusion-decreasing sequence
of ideals which does not stabilize. Since $\dot f$ is forced unbounded, $J=\bigcap J_n$. On the other hand, suppose that $J=\bigcap_nJ_n$ is the intersection
of an inclusion decreasing sequence of ideals which does not stabilize. Each ideal $J_n$ is dense in $Q_J$, so we can find a maximal antichain
$A_n\subset J_n$. Now suppose $a\in Q_J$. Either there is no condition $b\subset a$ which is compatible with at most countably many elements of $\bigcup_nA_n$.
In such a case $a\Vdash$``$\aleph_1$ is collapsed and by the CH an unbounded real is added''. Or there is such a condition $b$, compatible only with elements
$\{a_n^k:k\in\gw\}$ of the antichain $A_n$. Let $\dot f\in\baire$ be defined by $f(n)=$ the unique $k$ such that $a_n^k$ is in the generic filter.
The condition $b$ forces this function to be well-defined, and we will be done if we prove that it forces it not to be bounded by any ground model function.
Indeed, if $c\subset b$ is a condition and $g\in\baire$ is a function, it must be the case that $c\notin J_n$ for some number $n$, $d=\bigcup_{k\in g(n)} a_n^k\in J_n$,
$c\setminus d\notin J_n$ and clearly $c\setminus d\subset c$ is a condition forcing $g(n)\leq\dot f(n)$.

Finally, for the fourth item, suppose that $J=\bigcap_n J_n$ is an intersection of an inclusion decreasing sequence of P-ideals which does not stabilize.
Every ideal $J_n$ is dense in the factor $Q_J$, therefore we can find a maximal antichain $A_n\subset J_n$. We will be done if we show that every condition
$a\in Q_J$ is compatible with uncountably many elements of one of these antichains. Indeed, if $a\in Q_J$, then $a\notin J_n$ for some number $n$, and $a$ must 
be compatible with uncountably many  elements of the antichain $A_n$. This is true because if $X\subset A_n$ is a countable set then there is a set
$b\in J_n$ containing all elements of $X$ modulo finite, and then $a\setminus b\notin J_n$ is a condition which is incompatible with all elements of the set $X$!
\end{proof}

The following example answers a question of Ilijas Farah--Question 4.3 of \cite{farah:gaps}.

\begin{example}
An analytic P-ideal $J$ such that the factor $Q_J$ collapses $\aleph_1$. Let $\langle\ga_n:n\in\gw\rangle$ be a decreasing sequence of positive real numbers
smaller than $1$. Let $J_n$ be the summable P-ideal associated with the weight function $k\mapsto k^{-\ga_n}$. We claim
that $\langle J_n:n\in\gw\rangle$ is an inclusion-decreasing sequence of ideals which does not stabilize. The inclusions are clear. To see that
stabilization is impossible, choose a number $n$ and a set $a\notin J_n$. We will produce a set $b\subset a$, $b\in J_n\setminus J_{n+1}$. By induction
on $m\in\gw$ find mutually disjoint finite sets $b_m\subset a$ such that $\Sigma_{k\in b_m}k^{-\ga_n}\leq 2^{-m}$ while $\Sigma_{k\in b_m}k^{-\ga_{n+1}}\geq 1$.
Then $b=\bigcup_m b_m$ will be as desired. To find the set $b_m$, first find a number $k_m\in\gw$ such that for every $k>k_m$ it is the case that $k^{-\ga_n}\leq
2^{-m-1}k^{-\ga_{n+1}}\leq 2^{-m-1}$ and then find a finite set $b_m$ consisting of numbers larger than $k_m$ such that the sum $\Sigma_{k\in b_m}k^{-\ga_n}$
is between $2^{-m-1}$ and $2^{-m}$.

Now let $J=\bigcap_nJ_n$. This is an $F_{\gs\gd}$ ideal, and a simple diagonalization arguments show that it is a tall P-ideal. The Proposition shows that
the factor $Q_J$ collapses $\aleph_1$ in the presence of CH. If CH fails, the argument only shows that $Q_J$ is not proper, and we do not know if it has to collapse
$\fc$ to $\aleph_0$.

Note that this ideal is of minimal possible complexity for the factor $Q_J$ to be improper. All factors of $F_\gs$ forcings are $\gs$-closed by a theorem of Ilijas Farah.
\end{example}

\begin{example}
Let $\langle K_n:n\in\gw\rangle$ be a decreasing sequence of ideals on $\gw$ which does not stabilize. Consider the forcings $L_n$ of all trees $T\subset\wtree$
such that every node $s\in T$ longer than some fixed $t\in T$ splits into $K_n$-positively many immediate successors. It is not difficult to show that
the posets $L_n$ are proper and have the \crn--the arguments closely follow those for Laver forcing. Let $I_n:n\in\gw$ be the ideals on $\baire$ associated
with these forcings; a Borel set $B\subset\baire$ is $I_n$-positive if and only if $[T]\subset B$ for some tree $T\in L_n$. It is not difficult to see
that the ideals form an inclusion-decreasing sequence which does not stabilize. Let $J_n=\trace(I_n)$, let $I=\bigcap_nI_n$, and let $J=\bigcap_nJ_n=\trace(I)$.
Since $P_I$ regularly embeds into $Q_J$ by Claim~\ref{embeddingclaim} and $P_I$ is not proper by the Proposition, it must be the case that $Q_J$ is not proper either.
\end{example}

\section{The complexity of the trace ideals}

There are a number of questions concerning the structure of the factor posets $Q_J$ for various simply definable ideals $J$ on $\gw$. In order to
address them directly, we must analyze the complexity of the trace ideals.

\begin{proposition}
\label{luproposition}
Suppose that $I$ is a $\gs$-ideal on $\cantor$ such that the factor forcing $P_I$ is proper and has the \crn, and every analytic $I$-positive set 
has a Borel $I$-positive subset. The following are equivalent:

\begin{itemize}
\item $I$ is \pioneoneonsigmaoneone--\cite{kechris:book} 29.E, 35.9.
\item the trace ideal $\trace(I)$ is co-analytic.
\end{itemize}
\end{proposition}

\begin{proof}
The top to bottom direction follows immediately from the definitions. For the bottom to top direction suppose that $\trace(I)$
is co-analytic. Let $B\subset \cantor$ be an analytic set, with a tree $T$ such that $B=p[T]$. The proof will be complete if we show that
$B\notin I$ is equivalent to the following analytic statement: there is a set $b\subset\bintree$ decomposed into
antichains $b=\bigcup_na_n$ and a function $g:b\to T$ for $n\in\gw$ such that

\begin{itemize}
\item the antichain $a_{n+1}$ refines $a_n$
\item $g$ preserves extension and whenever $u\in a_n$ then $g(u)$ is a sequence of length $n$ whose first coordinates form an initial segment of $u$
\item $b\notin\trace(I)$.
\end{itemize}

First, if there are such objects $b$ and $g$, it is clear that the $I$-positive set $\pi(b)$ is a subset of $B$ and therefore $B\notin I$. On the other hand,
if $B\notin I$ then $B$ has a Borel $I$-positive subset $C$ by the assumptions, and by the properness
and the \crn\ of the poset $P_I$ there even is a Borel $I$-positive $G_\gd$ set $D\subset C$ and
a continuous function $f:D\to [T]$ such that for every real $r\in D$, the first coordinate of the value $f(r)$ is just $r$ itself.
It is then easy to construct $b$ and $g$ as above in such a way that $\pi(b)=D$ and for every real $r\in D$, $f(r)=\bigcup_{u\subset r} g(u)$.
\end{proof}

This lemma gives us a rather good criterion for checking whether a given trace ideal is co-analytic or not. If the poset $P_I$ adds a dominating real
then the ideal $I$ is not \pioneoneonsigmaoneone, \cite{z:book}, C.0.16. A quick review of forcings used in practice shows that many of them 
which do not add dominating reals (such as the Cohen or Solovay real) are associated with \pioneoneonsigmaoneone\ ideals.
However, Arnold Miller \cite{miller:snandb} constructed a definable c.c.c. ideal $I$ such that the poset $P_I$ does not add a dominating real while 
the ideal $I$ is still not \pioneoneonsigmaoneone. Therefore a careful check for this property is frequently necessary.

\begin{example}
\label{complexityexample}
Suppose that $J$ is an ideal on $\gw$ containing all finite sets, and let $I$ be the $\gs$-ideal on $\cantor$ associated with the
Prikry forcing $P(J)$ from Example~\ref{prikryexample}. The ideal $I$ is \pioneoneonsigmaoneone\ if and only if the ideal $J$ is $F_\gs$.
\end{example}

\begin{proof}
For the right-to-left direction note that the ideal $I$ is ergodic in the sense of \cite{z:book}: for every Borel set $B\subset\cantor$
closed under finite changes, either $B\in I$ or $\cantor\setminus B\in I$. Now if $I$ is \pioneoneonsigmaoneone\ then
the restriction of the trace $\trace(I)$ to trees on $\bintree$ is analytic: $T\in\trace(I)$ if and only if $[T]\in I$ if and only if $\cantor\setminus$ the closure of
the set $[T]$ under finite changes is $I$-positive, which is an analytic statement. Analytic ideals of compact sets are $G_\gd$ by a theorem
of Kechris, Louveau, and Woodin \cite{klw:compact}. Now it is not difficult to see that for a set $a\subset\gw$, $a\in J$ iff $T_a\notin\trace(I)$, where
$T_a$ is the set of all sequences $t\in\bintree$ for which $\forall n\ t(n)=1\liff n\notin a$. However, $T_a\notin\trace(I)$ is an $F_\gs$ condition!

For the left-to-right direction, write $\jx=\{a\subset\gwgw\setminus\{0\}: \exists b\in J\forall x\in a\ x\cap b\neq 0\}$.
It is clear that this is an ideal on the set $\gwgw\setminus\{0\}$. A useful observation:

\begin{claim}
If $J$ is an $F_\gs$ ideal then $\jx$ is $F_\gs$ again.
\end{claim}

\begin{proof}
By a theorem of Mazur \cite{mazur:fsigma} there is a lower semicontinuous submeasure $\mu$ on $\cP(\gw)$ such that $J=\{a\subset\gw:\mu(a)<\infty\}$. Let $\mu^{<\gw}$ be
a function on $\cP(\gwgw)$ defined by $\mu^{<\gw}(b)=\inf\{\mu(a):\forall x\in b\ x\cap a\neq 0\}$. It is not difficult to verify that
this is a lower semicontinuous submeasure such that $\jx=\{b\subset\gwgw:\mu^{<\gw}<\infty\}$. The claim follows.

In fact the proof of the Example shows that if the ideal $J$ is not $F_\gs$ then the ideal $\jx$ is not even analytic.
\end{proof}

By \cite{z:book} C.0.14, to prove the
left-to-right implication of the Example it is just necessary to show that the collection of countable subsets of $P(J)$ which are maximal antichains is a Borel set.
In order to do this, let $A\subset P(J)$ be a countable set. Then $A$ is a maximal antichain if and only if it is an antichain and for every finite set
$t\subset\gw$, every condition of the form $\langle t,a\rangle$ is compatible with some element of $A$. The latter condition is equivalent
to: either there is some condition $\langle u, b\rangle\in A$ such that $u\subset t$ and $b\cap t\setminus u=0$, or the set $a_t=\{x\subset\gw:
\exists b\ \langle t\cup x, b\rangle\in A\}$ is not in the ideal $\jx$. By the Claim, this is a Borel statement.
\end{proof}

On the other hand, we do not have a good criterion as to when the trace ideal is analytic. We have just a conjecture:

\begin{conjecture}
Suppose that $I$ is a $\gs$-ideal on $\cantor$ or $\baire$ such that the factor poset $P_I$ is proper with \crn. If the trace ideal
is analytic then it is in fact Borel.
\end{conjecture}

This conjecture can be viewed as a variation on the Kechris-Louveau-Woodin theorem \cite{klw:compact}. We can verify it in a good number of cases:

\begin{lemma}
Suppose that $I$ is a $\gs$-ideal on $\cantor$ such that $P_I$ is proper and bounding. If the trace ideal is analytic then it is Borel.
\end{lemma}

\begin{proof}
Since the poset $P_I$ is bounding,
compact sets are dense in it and it has the \crn--\cite{z:book}. 
If a set $a\subset\bintree$ is not in the trace ideal, apply these two properties below the condition
$\pi(a)$ to the name $\dot f(n)=$the $n$-th initial segment of the generic real in the set $\check a$. It follows that a set $a\subset\bintree$ is not in the trace ideal
if and only if there is a tree $T$ and disjoint finite subsets $\{b_n:n\in\gw\}$ of $T$
such that $[T]\notin I$ and each $b_n$ is a maximal antichain in $T$ consisting only of elements of the set $a$.
What is the complexity of the latter statement? The trace ideal restricted to trees is analytic, therefore
Borel by the Kechris-Louveau-Woodin theorem \cite{klw:compact}, and so this is an analytic statement. Thus the trace ideal is both analytic and co-analytic, therefore Borel.
\end{proof}

Of course, it is quite interesting to investigate the possibility of the trace ideal being Borel. The classical examples of such behavior are the meager and null ideals.
Other examples turn out to be quite hard to find.

\begin{example}
Suppose that $J$ is an analytic P-ideal, and $I$ is the $\gs$-ideal associated with the Prikry poset $P(J)$. Then the following are
equivalent:

\begin{enumerate}
\item $J$ is $F_\gs$
\item $I$ is \pioneoneonsigmaoneone
\item $P(J)$ does not add a dominating real
\item the trace ideal is Borel.
\end{enumerate}
\end{example}

\begin{proof}
(1) is equivalent with (2) by Example~\ref{complexityexample}. (2) implies (3) by \cite{z:book}, C.0.16. (3) implies (1) by a result of Solecki:
if an analytic P-ideal $J$ is not $F_\gs$ then the ideal $0\times$Fin is Rudin-Blass reducible to $J$
\cite{solecki:analyticideals}, 3.3. Let $f:\gw\to\gw\times\gw$ be such a finite-to-one reduction.
Since the $P(J)$ generic set $\dotagen\subset\gw$ has a finite intersection with a ground model set $b\subset\gw$ if and only if
$b\in J$, it immediately follows that $f''\dotagen$ has a finite intersection with a ground model set $b\subset\gw\times\gw$ if and only if
$b\in 0\times$Fin. Let $g\in\baire$ be defined by $g(n)=\min\{m\in\gw:\langle n,m\rangle\in f''\dotagen\}$. A brief inspection reveals
that this is a well-defined function modulo finite dominating all ground model functions.

This leaves us with the equivalence of (2) and (4). Note that the forcing $P_I$ has the \crn\ by Example~\ref{prikryexample} and so (4) implies (2)
by Proposition~\ref{luproposition}. For the opposite direction note that (2) implies the trace ideal is co-analytic by that same Proposition,
so it is enough to show from (2) that the trace ideal is analytic. 
To this end, use the ergodicity of the ideal $I$ again. For a set $a\subset\cantor$, $a\in\trace(I)$ if and only if the complement of the closure of the set
$\pi(a)$ under finite changes is $I$-positive, which is an analytic condition by (2).
\end{proof}

In order to produce many Borel ideals on $\gw$ such that their factors are mutually distinct as forcing notions, it is now necessary to show
that the posets $P(J)$ are distinct as forcings for different $F_\gs$ P-ideals $J$. This does not seem to be easy but we do have a natural criterion.

\begin{definition}
A forcing $P$ \emph{separates} an ideal $K$ on $\gw$ if it introduces a set $x\subset\gw$ which has finite intersection with all $K$-small sets
and infinite intersection with all $K$-positive sets.
\end{definition}

It is natural to suggest that the reason two forcings $P(K)$ and $P(J)$ may be different is because $P(J)$ does not separate $K$.
We will produce a forcing relation free criterion for this situation. First, a couple of basic facts.

\begin{proposition}
Suppose that $K,J$ are ideals on $\gw$.

\begin{enumerate}
\item If $K$ is Rudin-Keisler reducible to $J$ and a forcing $P$ separates $J$ then it separates $K$
\item $P(J)$ separates $J$, $P(J)\times$Cohen separates $\jx$
\item if $P$ separates $\jx$ then the poset $P(J)$ regularly embeds into $P\times$Cohen.
\end{enumerate}
\end{proposition}

\begin{proof}
Suppose $f\in\baire$ is a Rudin-Keisler reduction of $K$ to $J$, and a forcing $P$ separates $J$ by introducing an infinite set $a\subset\gw$.
It is straightforward to check that $f''a$ separates $K$.

For (2), let $a\subset\gw$ be a $P(J)$-generic set
and $f\in\baire$ be a Cohen real generic over $a$. A straightforward density argument shows that the set $a$ separates the ideal $J$. 
Another density argument will show that the set $b=\{x\subset a:x$ is nonempty, finite, and $\max(x)\in f(\min(x))\}$ separates $\jx$.

Finally, for (3) assume that $b\subset\gwgw$ is a set separating the ideal $\jx$ and $f\in 2^b$ is a Cohen real generic over $b$, obtained
by the poset $R$ of all finite partial functions from $b$ to $2$. We claim that the set $a=\bigcup\{x\in b:f(x)=1\}$ is $P(J)$-generic over the ground model.
For an ease of notation, if $p\in R$ is a condition, write $t_p=\bigcup\{x\in\dom(p):p(x)=1\}$.
Now suppose that $D\subset P(J)$ is a ground model open dense set and $p\in R$ is a condition. By
a standard genericity argument, it will be enough to find
an extension $q\supset p$ such that there is a $J$-small set $a\in V$ such that the condition $\langle t_q, a\rangle\in P(J)$
is in the open dense set $D$, and every set in $b\setminus\dom(q)$ has an empty intersection with the set $a$.

The set $c=\{x\in\gwgw:\exists a\in J\ \langle t_p\cup x, a\rangle\in D\}\in V$ is $\jx$-positive, therefore has infinite intersection with the set $b$.
Let $x\in b\cap c\setminus\dom(p)$ be an arbitrary set, and let $a\in J$ be a set such that $\langle t_p\cup x, a\rangle\in D$. Since the set $b$ separates
the ideal $\jx$, only finitely many of its elements have nonempty intersection with the set $a$. Let $q\supset p$ be a condition whose domain includes
all of these elements and $x$, and $x$ is the only element of $\dom(p\setminus q)$ on which the function $q$ returns value $1$. The condition $q\in R$ is as required.
\end{proof}

Here is the promised criterion:

\begin{proposition}
Suppose $K$ is a tall P-ideal and $J$ is an ideal on $\gw$. Then (1) implies (2) implies (3), where
\begin{enumerate}
\item $P(J)$ separates $K$
\item $K$ is Rudin-Keisler reducible to the ideal $\jx\restriction b$ where $b\subset\gwgw$ is some $\jx$-positive set not containing the empty set.
\item $P(J)\times$Cohen separates $K$
\end{enumerate}
\end{proposition}

\begin{proof}
(2) to (1) is proved in the contrapositive. 
Suppose (2) fails and $\dot x$ is a $P(J)$-name for an infinite set of natural numbers. We must show that it does not separate the ideal $K$.
Fix a partition $\gw=\bigcup_n z_n$ into infinite sets. For every set $t\in\gwgw$ and every number $n$ let $D(t,n)$ be the collection of those
sets $s\in\gwgw$ such that $s\neq 0, t\cap s=0$, $|s|\in z_n$ and there is a set $a\in J$ such that the condition $\langle t\cup s,a\rangle$ decides the
value of the $n$-th element of the set $\dot x$. Note that for each $t,n,s$ there is just one possible value because such conditions are mutually compatible; denote that
value by $f(t,n)(s)$. 

\begin{claim}
For every $t\in\gwgw$ and every $n\in\gw$ the set $D(n,t)$ is $\jx$-positive.
\end{claim}

By the assumptions, for every $t,n$ the function $g(t,n):\bigcup_{m>n}D(t,m)\to\gw$ defined by $g(t,n)=\bigcup_{m>n}f(t,m)$ is not a Rudin-Keisler reduction
of the ideal $K$. There are two cases. 

Either there is a set $t$ such that for some number $n$ there is a $K$-positive set $y\subset\gw$ such that $g(t,n)^{-1}y\in\jx$. In such a case, find
a set $a\in J$ which intersects every element of the preimage $g(t,n)^{-1}y$. It is not difficult to see that the condition
$\langle t,a\rangle$ forces all except the first $n$ elements of the set $\dot x$ to miss the set $y$, and so $\dot x$ does not separate the ideal $K$.

Or else for every set $t$ and for every number $n$ there is a $K$-small set $y(t,n)\subset\gw$ such that $g(t,n)^{-1}y(t,n)\notin\jx$. 
The sets $t\in\gwgw$ are now divided into two subcases. In the subcase A, if there
are infinitely many numbers $n$ such that the set $y(t,n)$ can be chosen to be a singleton then by the density of the ideal $K$ there is an infinite
set $y(t)\subset\gw$ in the ideal $K$ such that for every number $n$ there is an element $m\in y(t)$ such that $k>n$ and $g(t,0)^{-1}\{m\}\notin\jx$.
In the subcase B, if there are only finitely many numbers $n$ such that the set $y(t,n)$ can be chosen to be a singleton then use the fact that the ideal $K$ is a P-ideal
and find a set $y(t)\in K$ which modulo finite includes all the sets $y(t,n)$. Finally, find a set $y\in K$ which modulo finite includes all the sets
$y(t)$. The proof will be complete once we show that $P(J)$ forces the intersection $\dot x\cap\check y$ to be infinite.

Well, choose a condition $\langle t, b\rangle$ and a number $n$. We must find an extension of the condition which forces some number $m>n$ in the set $y$
into the set $\dot x$. First suppose that the set $t$ falls into the subcase A. In such a case, find a number $m>n$ in $y\cap y(t)$ such that 
$g(t,0)^{-1}\{m\}\notin\jx$. For such a number $m$ there is a set $s\in\gwgw$ disjoint from $b$ and a condition $\langle t\cup s, a\cup b\rangle\leq\langle t,b\rangle$
which forces $\check m\in\dot x$ and we are done. Now suppose that the set $t$ falls into the subcase B. Let $k>n$ be such that 
the set $y(t,k)$ had to be chosen infinite and $y(t,k)\setminus y$.
The $g(t,k)$-preimage of $y(t,k)\setminus(y\cup n)$ is then $\jx$-positive and there must be a set $s\in\gwgw$ disjoint from the set $b$ and a condition
$\langle t\cup s, b\cup a\rangle$ which forces some number $m\in y(t,k)\setminus n\subset y\setminus n$ into the set $\dot x$.

(2) implies (3) by the previous proposition. The forcing $P(J)\times$Cohen separates $\jx$, therefore it separates every ideal
$\jx\restriction b$ where $b$ is a positive set, and also every ideal Rudin-Keisler reducible to it.
\end{proof}

Perhaps some remarks are in order. The statements (1) and (3)
in the previous proposition seem to be very close but perhaps not equivalent.
We could not find and verify an example of an ideal $J$ such that the forcing $P(J)$ does not separate the ideal $\jx$. Such an ideal $J$ would have to be
very strongly inhomogeneous, and in fact any tall summable ideal is a plausible candidate. On the other hand, there are many ideals $J$
such that the forcing $P(J)$ does separate $\jx$. A good example of such a behavior is any ideal of the form $J=K^{<\gw}$ for some ideal $K$
or the ideal of nowhere dense subsets of the rationals.

\section{Open questions}

When the second author's stay in Mexico drew to a close, we decided to finish the work on this paper and leave several interesting questions open.

\begin{question}
Let $I$ be the meager ideal on $\cantor$, let $J=\trace(I)$. What is the remainder forcing $Q_J/P_I$? Similarly for the Lebesgue measure zero ideal.
\end{question}

\begin{question}
Is there a simple preservation criterion on the forcing $P_I$ which is equivalent to the remainder forcing being equal to $\pomodfin$?
\end{question}

\begin{question}
Is every proper forcing of the form $P_I$ regularly embeddable into a proper forcing of the form $Q_J$? Does every proper forcing of the form
$P_I$ have a presentation with the \crn?
\end{question}

\begin{question}
Assume CH. Is it true that for every ideal $J$ on a countable set, the factor forcing $Q_J$ collapses $\aleph_1$ if and only if it is $\aleph_0$-generated?
\end{question}

\begin{question}
The various definable improper forcings produced in the paper should collapse $\fc$ to $\aleph_0$ in ZFC. Is this really true?
\end{question} 

\begin{question}
For which $F_\gs$ P-ideals $J,K$ it is the case that $K$ is not Rudin-Keisler reducible to the ideal $\jx\restriction b$ for any $\jx$-positive set
$b\subset\gwgw$?
\end{question}

\bibliography{hz2}

\end{document}